\documentclass{amsart}
\usepackage{amssymb,  amscd, times}

\def\CC{{\mathbb C}}  
 
\def\FF{{\mathbb F}} 
 
\def\HH{{\mathbb H}}

\def\PP{{\mathbb P}}
\def\QQ{{\mathbb Q}}

\def\ZZ{{\mathbb Z}}

\def\G{\Gamma}
 
\def\hor{{\rm hor}}

\def\ver{{\rm vert}}

\newcommand{\p}{\partial}

\def\Acal{{\mathcal A}}

\def\Ccal{{\mathcal C}}

\def\Fcal{{\mathcal F}} 
\def\Hcal{{\mathcal H}} 
\def\Ical{{\mathcal I}}

\def\Lcal{{\mathcal L}}

\def\Ocal{{\mathcal O}}

\def\g{\gamma}
\def\la{\langle}
\def\ra{\rangle}

\def\cfrak{\mathfrak{c}}
\def\hfrak{\mathfrak{h}}
\def\gfrak{\mathfrak{g}}

\def\lfrak{\mathfrak{l}}
\def\lfrakhat{\hat{\mathfrak{l}}}

\def\Lg{L {\mathfrak{g}}}
\def\Lghat{\widehat{L\mathfrak{g}}}

\def\Lcalhatg{\widehat{\Lcal\mathfrak{g}}}

\def\afrak{\mathfrak{a}}
\def\mfrak{\mathfrak{m}}

\def\Ubar{{\overline{U}}}

\def\End{\operatorname{End}}
\def\pt{{\scriptscriptstyle\bullet}}

\def\half{{\tfrac{1}{2}}}

\newcommand\ad{\operatorname{ad}}

\newcommand\ann{\operatorname{ann}}

\newcommand\Gr{\operatorname{Gr}}

\newcommand\Hom{\operatorname{Hom}}

\newcommand\res{\operatorname{Res}}
\newcommand\spec{\operatorname{Spec}}

\newcommand\supp{\operatorname{Supp}}
\newcommand\sym{\operatorname{Sym}}

\newtheorem{theorem}{Theorem}[section]
\newtheorem{lemma}[theorem]{Lemma}
\newtheorem{proposition}[theorem]{Proposition}

\newtheorem{corollary}[theorem]{Corollary}

\theoremstyle{definition}

\theoremstyle{remark} 
\newtheorem{remark}[theorem]{Remark}

\title{Conformal blocks revisited}
\author{Eduard Looijenga}
\address{Mathematisch Instituut\\
Universiteit Utrecht\\
P.O.~Box 80.010, NL-3508 TA Utrecht\\
Nederland}
\email{looijeng@math.uu.nl}
\begin{document}

\begin{abstract}
We give a simple coordinate free description of the WZW connection and 
derive its main properties.
\end{abstract}
\date{July 4 2005}
\maketitle
This paper is based on a seminar talk whose original goal was to define 
the (projective) Wess-Zumino-Witten connection and to show its  flatness. My reason for writing this up was to demonstrate  the simplicity of its construction, a fact
not so apparent from the literature known to me, and to  recover  
in a selfcontained (and hopefully, more transparant) manner all the main results
of the fundamental paper by Tsuchiya-Ueno-Yamada \cite{tuy}.
Because of its consistent coordinate free approach, the present discussion 
is mostly algebraic in character, which (as so often) is not only conceptually  more satisfying, but also tends to simplify the arguments. (Another such approach insofar it involves the  WZW connection, can be found in a paper by Y.~Tsuchimoto \cite{tsuchimoto}.)

Sections 1 and 2 basically redo parts of the book of 
Kac-Raina \cite{kacraina} in a more algebraic setting,  replacing for 
instance the ring of complex Laurent polynomials by a complete local
field containing $\QQ$ (or rather a direct sum of these), which 
is then also allowed to `depend on parameters'. What may be new is the  
extension \ref{cor:derive} of the Sugawara representation to a relative 
situation involving a Leibniz rule in 
the horizontal direction. I regard this construction as the origin
of WZW-connection and its projective flatness. The connection
is defined in Section 3. This is followed by a derivation of the coherence 
of the sheaf of conformal blocks and 
what is called the propagation of vacua.   
Special attention is paid to the genus zero case and it 
shown how the WZW-connection is then related to the one of 
Knizhnik-Zamolodchikov. Finally, we derive in 
Section 4 the basic results associated to a double point 
degeneration (such as local freeness, factorization and monodromy).

The reader be warned that the bibliography is merely a record of 
papers consulted, rather than anything else. 

\section{Canonical construction of the Virasoro algebra}\label{section:vir} 
In this section we fix a $\QQ$-algebra $R$ and
a $R$-algebra $\Ocal$ isomorphic to the formal power series ring $R[[t]]$.
In other words, $\Ocal$  comes with a principal ideal $\mfrak$
so that $\Ocal$ is complete for the $\mfrak$-adic topology and $\Ocal/\mfrak^j$ is
for $j=1,2,\dots$ a free $R$-module of rank $j$. A generator $t$ of $\mfrak$
will then identify $\Ocal$ with $R[[t]]$. We denote by
$L$ the localization of $\Ocal$ obtained by inverting a generator of $\mfrak$.
For $N\in\ZZ$, $\mfrak^N\subset L$ has the obvious meaning
as a $\Ocal$-submodule of $L$. The \emph{$\mfrak$-adic topology} on $L$ is the topology that has the collection of cosets $\{ f+\mfrak^N\}_{f\in L,N\in
\ZZ}$ as a basis of open subsets. We sometimes write $F^NL$ for $\mfrak^{N}$. We further denote by
$\theta$ the $L$-module of continuous $R$-derivations from $L$ into $L$
and by $\omega$ the $L$-dual of $\theta$. These $L$-modules 
come with filtrations (making them principal filtered $L$-modules): $F_N\theta$ consists of the derivations that take $\mfrak$ to 
$\mfrak^{N+1}$ and $F^N\omega$ consists of the $L$-homomorphisms $\theta\to L$ that take  $F^0\theta$ to $\mfrak^N$.
So in terms of the generator $t$ above, $L=R((t))$, $\theta=R((t))\frac{d}{dt}$, 
$F^N\theta=R[[t]]t^{N+1}\frac{d}{dt}$,  $\omega =R((t))dt$ and $F^N\omega=R[[t]]t^{N-1}dt$. 

The residue map
$\res: \omega\to R$ which assigns to an element of $R((t))dt$ the coefficient of $t^{-1}dt$
is independent of the choice of $t$. The $R$-bilinear map
\[
r: L\times \omega \to R,\quad  (f,\alpha)\mapsto \res (f\alpha)
\]
is a topologically perfect pairing of filtered $R$-modules: we have $r(t^k, t^{-l-1}dt)=\delta_{k,l}$ so that any $R$-linear  $\phi: L\to R$  which is continuous (i.e., $\phi$ zero on $\mfrak^N$ for some $N$) is definable by an element of $\omega$
(namely by $\sum_{k<N} \phi (t^k)t^{-k-1}dt$) and likewise
for a $R$-linear continuous map $\omega\to R$.

\subsection*{A trivial Lie algebra}
If we regard $L^\times$ as an algebraic group over $R$ (or rather as a group object
in a category of ind schemes over $R$), then
its Lie algebra, denoted here by $\lfrak$, is $L$,  regarded as a $R$-module with trivial Lie bracket. It comes with a decreasing filtration $F^\pt\lfrak$ (as a Lie algebra) defined by the valuation.
The universal enveloping algebra $U\lfrak$ is clearly 
$\sym_R^\pt (\lfrak)$. 
The ideal $U_+\lfrak\subset U\lfrak$ generated by $\lfrak$
is also a right  $\Ocal$-module (since $\lfrak$ is). We complete it $\mfrak$-adically:
given an integer $N\ge 0$, then an $R$-basis of 
$U_+\lfrak/(U\lfrak\circ F^N\lfrak)$ 
is the collection $t^{k_1}\circ\cdots\circ t^{k_r}$ with 
$k_1\le k_2\le \cdots \le k_r<N$. So elements of the  completion 
\[
U_+\lfrak\to \Ubar_+\lfrak:=\varprojlim_N U_+\lfrak/U\lfrak\circ F^N\lfrak
\]
are series of the form 
$\sum_{i=1}^\infty r_i t^{k_{i,1}}\circ\cdots\circ t^{k_{i,r_i}}$
with $r_i\in R$, $k_{1,i}\le k_{2,i}\le\cdots \le k_{i,r_i}$, $\{ k_{i,1}\}_i$ bounded from below and $\lim_{i\to\infty}k_{i,r_i}=\infty$. 
We put $\Ubar\lfrak:= R\oplus \Ubar_+\lfrak$, which could  of course 
have been defined directly as
 \[
U\lfrak\to \Ubar\lfrak:=\varprojlim_N U\lfrak/U\lfrak\circ F^N\lfrak.
\]
We will refer to this construction as the \emph{$\mfrak$-adic completion on the right}.
(In the present case, this is no different from  $\mfrak$-adic completion 
on the left, because 
$\lfrak$ is commutative.)

Any  $D\in\theta$ defines an $R$-linear map $\omega\to L$ which is 
selfadjoint relative our topological pairing: $r(\la D,\alpha \ra,\beta)=r(\alpha,\la D,\beta\ra)$.
We use that pairing to identify $D$ with an
element of the closure of $\sym^2\lfrak$ in $\Ubar\lfrak$. Let $C(D)$ be half this element, so that in terms of the above topological basis:
\[
C(D)=\half \sum_{l\in\ZZ} \la D, t^{-l-1}dt\ra\circ t^l.
\]
 In particular we have for $D=D_k=t^{k+1}\frac{d}{dt}$
\[
C(D_k)=\half \sum_{i+j=k} : t^i\circ t^j :\quad.
\]
We here adhered to the normal ordering convention (the factor with the 
highest index comes last and hence acts first), to make the righthand 
side look like an element of $\Ubar\lfrak$, although
there is no need for this as  $t^i\circ t^j=t^j\circ t^i $. Observe 
that the map  $C: \theta\to\Ubar\lfrak$ is continuous. 

\subsection*{Oscillator algebra}
The residue map defines a central extension of $\lfrak$, the 
\emph{oscillator algebra} $\lfrakhat$, which as a $R$-module is 
simply $\lfrak\oplus \hbar R$ and has Lie bracket
\[
[f +\hbar r ,g+\hbar s ]:=\hbar\res (g\, df).
\]
So $[t^{k},t^{-l}]=\hbar k\delta_{k,l}$ and $R\hbar$ is central. We filter 
$\lfrakhat$ by letting $F^N\lfrakhat$ be $F^N\lfrak$ for $N>0$
and $F^N\lfrak +R\hbar$ for $N\le 0$. We complete 
$U\lfrakhat$ $\mfrak$-adically on the right:
\[
U\lfrakhat\to \Ubar\lfrakhat:=\varprojlim_N U\lfrakhat /U\lfrakhat\circ F^N\lfrak.
\]
Since $\hbar$ is in the center of $\lfrak$, this is a $R[\hbar ]$-algebra
(for a similar reason it is even a $R[e,\hbar ]$-algebra if $e=t^0$ denotes the unit element of $L$ viewed as an element of $\lfrak$).  As a $R[\hbar]$-algebra it is a quotient of the tensor algebra of 
$\lfrak$ (over $R$) tensored with $R[\hbar]$, $\otimes^\pt_R \lfrak \otimes_R R[\hbar]$,
 by the two-sided ideal generated by the elements $f\otimes g-g\otimes f-\hbar \res (gdf)$.  As a $R[\hbar]$-module it has for topological basis  the collection 
 $t^{k_1}\circ\cdots\circ t^{k_r}$ with $r\ge 0$, $k_1\le k_2\le \cdots \le k_r$. 
 Since $\hat\lfrak$ is not abelian, the left and  right $\mfrak$-adic topologies differ. For instance, 
$\sum_{k\ge 1} t^k\circ t^{-k}$ does not converge in $\Ubar\lfrakhat$, whereas 
$\sum_{k\ge 1} t^{-k}\circ t^k$ does. Notice that the obvious surjection 
$\pi: U\lfrakhat\to U\lfrak$ is simply the reduction modulo $\hbar$ of 
$U\lfrakhat$ and likewise for their completions. 
Observe that the filtrations of $\lfrak$ and $\hat\lfrak$ determine decreasing 
filtrations
of the their (completed) universal enveloping algebras. For instance,
$F^NU\hat\lfrak=\sum_{r\ge 0}\sum_{n_1+\cdots +n_r\ge N} 
F^{n_1}\hat\lfrak\circ\cdots \circ F^{n_r}\hat\lfrak$.

Denote by $\lfrak_2$ the image of 
$\lfrak^{\otimes 2}\subset\lfrakhat^{\otimes 2}\to U\lfrakhat$.
Under the reduction modulo $\hbar$,  $\lfrak_2$ maps onto $\sym^2(\lfrak)\subset U\lfrak$ with kernel $R\hbar$.
Its  closure $\bar \lfrak_2$ in $\Ubar\lfrakhat$ maps onto the closure of 
$\sym^2(\lfrak)$ in $\Ubar \lfrak$ with the same kernel. We denote  
by $\hat\theta$ the set of  pairs $(D, u)\in\theta\times \bar \lfrak_2$ 
for which $C(D)$ is the mod $\hbar$ reduction of $u$ so that we have an exact sequence
\[
0\to \hbar R \to \hat\theta \to \theta \to 0
\]
of $R$-modules and a natural $R$-homomorphism 
$\hat C:\hat\theta\to \Ubar\lfrakhat$. 
The generator $t$ defines a (noncanonical) section of $\hat\theta \to \theta$: 
\[
D\in\theta\mapsto \hat D:= (D,\half\sum_{j\in\ZZ} : \la D, t^{-j-1}dt\ra\circ t^j:) \in \theta\times\lfrakhat .
\]

\begin{lemma}\label{lemma:chat}
We have
\begin{enumerate}
\item[(i)] $[\hat C(\hat D),f]= -\hbar D(f)$ as an identity  in $\Ubar \lfrakhat$ 
(where $f\in\lfrak\subset\hat\lfrak$) and
\item[(ii)] $[\hat C(\hat D_k), \hat C(\hat D_l)]=-\hbar(l-k)  
\hat C(\hat D_{k+l})+\hbar^2\frac{1}{12}(k^3-k)\delta_{k+l,0}$.
\end{enumerate}
\end{lemma}
\begin{proof}
For the first statement we compute $[\hat C(\hat D_k), t^l]$. The only terms
in the expansion of $\sum_{i+j=k} : t^i\circ t^j:$ that can contribute 
are for the form $[t^{k+l}\circ t^{-l}, t^l]$ or $[t^{-l}\circ t^{k+l}, t^l]$,
(depending on whether $k+2l\le 0$ or $k+2l\ge 0$) and with coefficient 
$\half$ if $k+2l= 0$ and $1$ otherwise. In all cases the result is 
$\hbar lt^{k+l}= \hbar D_k(t^l)$.
 
Formula (i) implies that 
\begin{multline*}
[\hat C(\hat D_k), \hat C(\hat D_l)]=
\lim_{N\to\infty}\sum_{|i|\le N} \half\left( D_k(t^i)\circ t^{l-i}+ t^i\circ D_k(t^{l-i})\right)\\
=-\hbar \lim_{N\to\infty}\sum_{|i|\le N} 
\left( it^{k+i}\circ t^{l-i}+ t^i\circ (l-i)t^{k+l-i}\right)
\end{multline*}
This is up to a reordering equal to $ -\hbar (l-k)\hat C(\hat D_{k+l})$.
The terms which do not commute and are in the wrong order are those
for which $0<k+i=-(l-i)$ (with coefficient $i$) and for which
$0<i=-(k+l-i)$ (with coefficient $(l-i)$). This accounts for
the extra term $\hbar^2\frac{1}{12}(k^3-k)\delta_{k+l,0}$
\end{proof}

This lemma suggests we rescale $\hat C$ as
\[
T:=-\tfrac{1}{\hbar}\hat C: \hat\theta\to \Ubar\lfrakhat [\tfrac{1}{\hbar}]
\]
and write $c_0$ for  $(0,-\hbar)\in \hat\theta$ since then 
\begin{enumerate}
\item[(i)] $T$ is injective and maps $\hat\theta$ onto a Lie subalgebra of 
$\Ubar\lfrakhat[\frac{1}{\hbar}]$ with $c_0\in\hat\theta\mapsto 1$,
\item[(ii)] if we transfer the Lie bracket to $\hat\theta$, we find that
\[
[\hat D_k, \hat D_l]= (l-k)\hat D_{k+l}+\frac{k^3-k}{12}\delta_{k+l,0}c_0, 
\]
\item[(iii)]and $\ad_{T(\hat D)}$ leaves $\lfrak$ 
invariant (as a subspace of $\Ubar\hat\lfrak$) and acts on that subspace by derivation with respect to $D$, 
\end{enumerate}
Thus we get a central extension of $\theta$ with a one-dimensional 
center canonically isomorphic to $R$ ($c_0$ corresponds to $1\in R$), 
the \emph{Virasoro algebra} (of the $R$-algebra $L$).

\begin{remark}
An alternative coordinate free definition of the Virasoro algebra, based on the algebra of
pseudodifferential operators on $L$, can be found in \cite{bms}.
\end{remark} 

\subsection*{Fock representation}
It is clear that $F^\lfrak=F^1\lfrakhat$ is an abelian
subalgebra of $\lfrakhat$. We put
\[
\FF:= (U\lfrakhat /U\lfrakhat\circ F^1\lfrak) [\tfrac{1}{\hbar}].
\]
This is a representation of $\lfrak$ over $R$; at the same time it is a 
$R[e,\hbar,\hbar^{-1}]$-module and as such it is free with basis 
the collection $t^{-k_r}\circ\cdots\circ t^{-k_1}v_o$, where $v_o$ denotes
the image of $1$, $r\ge 0$ and $1\le k_1\le k_2\le \cdots\le k_r$ (for $r=0$, read $v_o$). This shows that we may identify
$\FF$ with $(\Ubar\lfrakhat/\Ubar\lfrakhat\circ F^1\lfrakhat)[\frac{1}{\hbar}]$, 
which makes it a representation for $\hat\theta$ over $R$. This is the 
\emph{Fock representation} of $\hat\theta$.  The 
$R[e,\hbar,\hbar^{-1}]$-subalgebra of $U\lfrakhat$ generated by the elements $t^{-k}$, $k=1,2,\dots$ is a polynomial algebra 
(in infinitely many variables) which projects isomorphically onto $\FF$. 
But this subalgebra depends on $t$; $\FF$ itself does not seem to be an 
$R$-algebra in a canonical way. On the other hand, it is naturally filtered
with $F^N\FF=F^N\Ubar\lfrakhat (v_0)$. 

It follows from Lemma \ref{lemma:chat} that when $D\in\theta$, then
\begin{multline*}
T(\hat D)t^{-k_r}\circ\cdots \circ t^{-k_1} v_0=\\
=\left(\sum_{i=1}^r X_{\alpha_r}t^{-k_r}\circ\cdots
D(t^{-k_i})\circ\cdots \circ t^{-k_1}\right) v_0+t^{-k_r}\circ\cdots \circ t^{-k_1} 
T_\gfrak(\hat D) v_0.
\end{multline*}
For $D\in F^0\theta$ we have $T_\gfrak(\hat D)v_0=0$ and so
$F^0\theta$ acts on $\FF$
by coefficientwise derivation.

\section{The Sugawara construction}\label{section:sugawara}

In this section, we fix a base field $k$ of caracteristic zero and
a simple Lie algebra $\gfrak$ over $k$ of finite dimension. We retain the data and the notation of Section \ref{section:vir}, except that
we now assume $R$ to be a $k$-algebra.

\subsection*{Loop algebras} The space of symmetric
invariant bilinear forms $\gfrak\times\gfrak\to k$ is of dimension
one and has a canonical generator $(\; |\; )$ (namely the form which takes the value
$2$ on the small coroots). We choose an orthonormal basis $\{ X_\kappa\}_\kappa$ 
of $\gfrak$ relative to this form.
The dual of this line is the space $\gfrak$-invariants in $\sym^2\gfrak$; 
we shall denote this line by $\cfrak$. The form $(\; |\; )$
singles out a generator of $\cfrak$, namely the Casimir element $c=\sum_\kappa X_\kappa\otimes X_\kappa$. It is well-known and easy to prove
that $\cfrak$ maps to the center of $U\gfrak$. In particular,
$c$ acts in any finite dimensional irreducible representation  of $\gfrak$ by a scalar.
In the case of the adjoint representation we denote half this scalar 
by $\check{h}$ (for this happens to be the \emph{dual Coxeter number} of 
$\gfrak$) so that
\[
\sum_\kappa [X_\kappa [X_\kappa, Y]]=2\check{h}Y\quad \text{for all $Y\in\gfrak$.}
\]
Let $L\gfrak$ stand for $\gfrak\otimes_k L$, but considered as a filtered $R$-Lie algebra (so we restrict the scalars to $R$):
$F^N L\gfrak= \gfrak\otimes_k \mfrak^N$.
An argument similar as for $r$ shows that the pairing
\[
r_\gfrak : (\gfrak\otimes_k L)\otimes (\gfrak\otimes_k\omega)\to R,\quad 
(Xf,Y\kappa)\mapsto (X|Y)\res (f\kappa)
\]
is topologically perfect; the basis dual to $(X_\kappa t^l)_{\kappa,l}$ is
$(X_\kappa t^{-l-1}dt)_{\kappa,l}$. 

For an integer $N\ge 0$, the
quotient $U_+L\gfrak/U L\gfrak\circ F^NL\gfrak$ is a free $R$-module with generators $X_{\kappa_1}t^{k_1}\circ\cdots\circ X_{\kappa_r}t^{k_r}$, $k_1\le\cdots\le k_r<N$. We complete $U L\gfrak$  $\mfrak$-adically on the right: 
\[
\Ubar L\gfrak:= \varprojlim_N  U L\gfrak/
U L\gfrak\circ  F^NL\gfrak.
\]
A central extension $\Lghat$ of $\Lg$  by $\cfrak$ is defined
by endowing $L\gfrak\oplus \cfrak$ with the Lie bracket
\[
[Xf+c r ,Yg+c s ]:=[X,Y]fg +c\res (g df) (X|Y).
\]
We filter $\Lghat$ by letting  for $N>0$,
$F^N\Lghat=F^N\Lg$ and for $N\le 0$, $F^N\Lghat=F^N\Lg +\cfrak$.
Then $U \Lghat$ is a filtered $R[c]$-algebra whose reduction modulo $c$ is 
$U L \gfrak$. Since the residue is zero on $\Ocal$, the inclusion of 
$F^0\Lg$ in $\Lghat$ is a homomorphism of Lie algebras.
The $\mfrak$-adic completion on the right
\[
\Ubar \Lghat:= \varprojlim_N  
U \Lghat/(U \Lghat\circ F^N\Lg)
\]
is still a $R[c]$-algebra and the obvious surjection 
$\Ubar\Lghat\to \Ubar L\gfrak$ is the reduction modulo $c$.
These (completed) enveloping algebras inherit a decreasing  filtration from $L$. 

\subsection*{Sugawara representation}
If $c$ is regarded as an element of $\gfrak\otimes_k\gfrak$, then tensoring with it defines the $R$-linear map
\[
\lfrak\otimes_R\lfrak\to \Lg\otimes_R \Lg,\quad
f\otimes g\mapsto \sum_{\kappa}X_\kappa f\otimes X_{\kappa}g,
\]
which, composed with $\Lg\otimes_R \Lg\subset \Lghat\otimes_R \Lghat\to U\Lghat$, yields $\g: \lfrak\otimes_R\lfrak\to U\Lghat$. Since 
$\g (f\otimes g-g\otimes f)=\sum_\kappa [X_\kappa f,X_\kappa g]=c\dim\gfrak \res (g df)$,  $\g$  drops and extends naturally to an $R$-module homomorphism  
$\hat\gamma: \lfrak_2\to U\Lghat $ which sends $\hbar$ to  $c\dim\gfrak$.
It extends continuously to a map from the closure $\bar\lfrak_2$ of
$\lfrak_2$ in $\Ubar \hat\lfrak$ to $\Ubar\Lghat$.
As $\bar\lfrak_2$ contains the image of $C:\hat\theta\to \Ubar\hat\lfrak$, we get a natural $R$-homomorphism 
\[
\hat C_\gfrak:=\hat\gamma\hat C :\hat\theta\to \Ubar \Lghat.
\]
We may  also describe $\hat C_\gfrak$ in the spirit of Section 
\ref{section:vir}: given $D\in \theta$, then the $R$-linear  map
\[
1\otimes D:  \gfrak\otimes_k\omega \to  \gfrak\otimes_k L
\]
is continuous and selfadjoint relative to $r_\gfrak$ and the perfect pairing 
$r_\gfrak$ allows us to identify it with an element of $\Ubar L\gfrak$; this element
is our $\hat C_\gfrak(D)$. 
The choice of $t$ yields the lift 
\[
\hat C_\gfrak(\hat D)=\half\sum_{\kappa}\sum_l  : 
X_\kappa\la D,t^{-l-1}dt\ra \circ X_\kappa t^l:\; \in \Ubar \Lghat 
\]
(so that in particular, 
$\hat C_\gfrak(\hat D_k)=\half\sum_{\kappa, l}: X_\kappa t^{k-l}\circ 
X_\kappa t^l:$). This formula can be used to define  $\hat C_\gfrak$, but this approach does  not exhibit its naturality.

\begin{lemma}
We have
\begin{enumerate}
\item[(i)]  $[C_\gfrak(\hat D), Xf]= -(c+\check{h} )XD(f)$, for every $D\in\theta$,
$X\in\gfrak$ and $f\in L$ and
\item[(ii)] $[C_\gfrak(\hat D_k), C_\gfrak(\hat D_l)]=(c+\check{h})
(k-l)C_\gfrak(\hat D_{k+l})+c(c+\check{h})\delta_{k+l,0} \frac{1}{12}\dim \gfrak (k^3-k)$.
\end{enumerate}
\end{lemma}

For the proof (which  is  a bit 
tricky, but not very deep), we refer to Lecture 10 of 
\cite{kacraina} (our  $C_\gfrak(\hat D_k)$ is their $T_k$). 

\begin{corollary}[Sugawara representation]\label{cor:sugawara}
The $R$-linear map 
\[
T_\gfrak:=\frac{-1}{c+\check{h}}\hat C_\gfrak: \hat\theta\to \Ubar\Lghat [ \tfrac{1}{c+\check{h}}]
\]
which  sends the central $c_0\in \hat\theta$ to $c(c+\check{h})^{-1}\dim \gfrak$ is a homomorphism
of $R$-Lie algebras. Moreover, if $D\in\theta$, then 
$\ad_{T_\gfrak(\hat D)}$ leaves $\Lg$ 
invariant (as a subspace of $\Ubar\Lghat$) and acts on that subspace by derivation with respect to $D$.
\end{corollary}

\subsection*{Fock type representation for $\gfrak$} 
Consider the $U\Lghat$-module
\[
\FF(L\gfrak):=(U\Lghat /U F^1\Lg)[\frac{1}{c+\check{h}}]
\]
If $v_0\in \FF(L\gfrak)$ denotes the image of $1$, then as  a $R[\frac{1}{c+\check{h}}]$-module it has for basis is the collection 
$\{ X_{\kappa_r}t^{-k_r}\circ\cdots \circ X_{\kappa_1}t^{-k_1} (v_0)\, :\,
r\ge 0,\quad  0\le k_1\le k_2\le \cdots \le k_r\}$.
It is easy to see that we can identify
$\FF(L\gfrak)$ with $\Ubar\Lghat /\Ubar F^1\Lg[\frac{1}{c+\check{h}}]$, 
so that $\hat\theta$ acts on $\FF(L\gfrak)$.
It follows from Corollary \ref{cor:sugawara} that when $D\in\theta$, then
\begin{multline*}
T_\gfrak(\hat D)X_{\kappa_r}t^{-k_r}\circ\cdots \circ X_{\kappa_1}t^{-k_1} (v_0)=\\
=\left(\sum_{i=1}^r X_{\kappa_r}t^{-k_r}\circ\cdots
X_{\kappa_i}D(t^{-k_i})\circ\cdots \circ X_{\kappa_1}t^{-k_1}\right) v_0+\\
+X_{\kappa_r}t^{-k_r}\circ\cdots \circ X_{\kappa_1}t^{-k_1} T_\gfrak(\hat D) v_0.
\end{multline*}
Thus $\hat\theta$ is  faithfully represented as a Lie algebra 
of $R[\frac{1}{c+\check{h}}]$-linear endomorphisms of $\FF(L\gfrak)$.
If $D\in F^0\theta$, then clearly $T_\gfrak(\hat D)v_0=0$ and hence:

\begin{lemma}\label{lemma:derive}
The Lie subalgebra $F^0\theta$ of $\hat\theta$ acts on $\FF(L\gfrak)$
by coefficientwise derivation.
\end{lemma}

This lemma has an interesting corollary. In order to state it, 
consider the module of $k$-derivations $R\to R$ (denoted here simply 
by $\theta_{R}$ instead of the more correct $\theta_{R/k}$)
and the module of $k$-derivations of $L$ which preserve $R\subset L$ 
(denoted by $\theta_{L,R}$). 
Since $L\cong R((t))$ as an $R$-algebra, every $k$-derivation 
$R\to R$ extends to  one from $L$ to $L$. 
So $\theta_{L,R}/\theta$ can be identified with $\theta_R$. 

\begin{corollary}\label{cor:derive}
The central extension $\hat\theta$ of  $\theta$ by $Rc_0$ naturally extends  to a
central extension of $R$-Lie algebras $\hat\theta_{L,R}$ of $\theta_{L,R} $ by $Rc_0$ 
(so with  $\hat\theta_{L,R}/\hat\theta =\theta_{L,R}/\theta\cong\theta_R$)
and the Sugawara representation $T_\gfrak$ of $\hat\theta$ on  $\FF(L\gfrak)$ 
extends to $\hat\theta_{L,R}$ in such a manner that it preserves any 
$U\Lghat$-submodule of $\FF(L\gfrak)$ and for every 
$D\in\theta_{L,R}$, the following relations hold in $\End (\FF(L\gfrak))$:
\begin{enumerate}
\item[(i)] $[T_{\gfrak}(\hat D),Xf]=X(Df)$ for $X\in\gfrak$, $f\in L$, and
\item[(ii)] (Leibniz rule) $[T_{\gfrak}(\hat D),r]=Dr$ for $r\in R$.
\end{enumerate}
These assertions also hold for the oscillator representation $T$ on $\FF$. 
\end{corollary}
\begin{proof}
The generator $t$ can be used to define a section of $\theta_{L,R}\to \theta_R$:
the set of elements of $\theta_{L,R}$ which kill $t$ is a $k$-Lie subalgebra
of  $\theta_{L,R}$ which projects isomorphically onto $\theta_R$.
Now if $D\in \theta_{L,R}$, write $D=D_\ver+D_\hor$ with $D_\ver\in \theta$ and 
$D_\hor (t)=0$ and define a $k[\frac{1}{c+\check{h}}]$-linear operator
$\hat D$ in $\FF(L\gfrak)$ as the sum of $T_\gfrak(\hat D_\ver)$ 
and  coefficientwise derivation by $D_\hor$. This map clearly satisfies the two 
properties and preserves any $U\Lghat$-submodule of $\FF(L\gfrak)$.
As to its  dependence on $t$: another choice yields a decomposition of the form $D=(D_\hor+D_0) +(D_\ver -D_0)$ 
with $D_0\in F^0\theta$ and in view of Lemma  
\ref{lemma:derive} $\hat D_0$ acts in $\FF(L\gfrak)$ by coefficientwise derivation. 

The same argument works for the oscillator representation.
 \end{proof}

\subsection*{Semilocal case}
This refers to the situation where we allow the $R$-algebra $L$ to be a finite direct sum of $R$-algebras isomorphic 
to $R((t))$: $L=\oplus_{i\in I} L_i$, where $I$ is a nonempty finite index set and $L_i$ as before. We use the obvious
convention, for instance, $\Ocal$, $\mfrak$, $\omega$, $\lfrak$ are now the direct sums over $I$ (as filtered objects)  of the items suggested by the notation.
If $r: L\times\omega\to R$ denotes the sum of the residue pairings of 
the summands, then $r$ is still topologically perfect. In this setting, 
the oscillator algebra $\hat\lfrak$ is of course not 
the direct sum of the $\hat\lfrak_i$, but rather
the quotient of  $\oplus_i\hat\lfrak_i$ that identifies the central generators 
$c_{0,i}$ of the summands with a single $c_0$. We thus get a Virasoro extension
$\hat\theta$ of $\theta$ by $c_0 R$ and a (faithful) oscillator representation of
$\hat\theta$ in $\Ubar\hat\lfrak$. The decreasing filtrations are the obvious
ones. In likewise manner we define $\Lghat$ (a central extension of 
$\oplus_{i\in I} L\gfrak_i$) and construct 
$\hat\theta_{L,R}$ and the associated Sugawara representation: 
Corollaries \ref{cor:sugawara} and  \ref{cor:derive} continue to hold.

\section{The WZW connection}\label{section:wzw}

From now on we assume that our base field $k$ is algebraically closed
and $R$ is a noetherian $k$-algebra.
We place ourselves in the semi-local case. 
In what follows, we often allow $\gfrak$ to be replaced by $k$, viewed as
an abelian Lie algebra, where the substitutions are the obvious ones;
for instance,  $\Lghat$,  $T_\gfrak$ and $\FF(L\gfrak)$ become $\lfrakhat$, $T$ and $\FF$. 

\subsection*{Abstract conformal blocks}
Let  $A$  be a $R$-subalgebra of $L$ and let $\theta_{A/R}$ have the usual meaning
as the Lie algebra of $R$-derivations $A\to A$. We assume that:
\begin{enumerate}
\item[($A_1$)] as a $R$-algebra, $A$ is flat and of finite type,
\item[($A_2$)] $A\cap \Ocal =R$ and  the $R$-module $L/(A+\Ocal)$ is locally 
free of finite rank,
\item[($A_3$)] the annihilator $\ann(A)$ of
$A$ in $\omega$ contains the image of $dA$ and the
resulting map $\Hom_A(\ann (A),A)\to\theta_{A/R}$ is an isomorphism.
\end{enumerate}

We denote $\theta_{A,R}$ the Lie algebra of $k$-derivations $A\to A$ 
which preserve $R$. So the quotient  $\theta_R^A:=\theta_{A,R}/\theta_{A/R}$ 
consists of the $k$-derivations  $R\to R$ that extend to one of $A $ 
and hence is a submodule of $\theta_R$. It is clear that 
$\theta_{A,R}\subset \theta_{L,R}$. We denote by $\hat\theta_{A,R}$ 
the preimage of  $\theta_{A,R}$ in $\hat\theta_{L,R}$ and by
$\hat\theta_R^A$ the quotient $\hat\theta_{A,R}/\theta_{A/R}$. 
These are central extensions of $\theta_{A,R}$ resp.\ $\theta_R^A$ 
by $\cfrak\otimes_k R =cR$.

We put 
\[
A\gfrak:= \gfrak\otimes_k A
\]
and view  it as a Lie subalgebra of $\Lg$.  Since the residue vanishes on 
$\ann(A)$, the inclusion $A\gfrak\subset \Lghat$ is a Lie algebra homomorphism. 
We define the \emph{universal covacuum space} as the space of 
$A\gfrak$-covariants in $\FF(L\gfrak)$, 
$\FF(L\gfrak)_{A\gfrak}:=\FF(L\gfrak)/A\gfrak \FF(L\gfrak)$.

\begin{proposition}\label{prop:A}
For  $D\in \theta_{A/R}$, $T_\gfrak(\hat D)$ lies in the closure of 
$A\gfrak\circ\Lghat$ in $\Ubar\Lghat$.  

The Sugawara representation of the Lie algebra $\hat\theta_{A,R}$ on 
$\FF(L\gfrak) $ preserves  the space of 
$A\gfrak$-covariants in $\FF(L\gfrak)$, $\FF(L\gfrak)_{A\gfrak}:=
\FF(L\gfrak)/A\gfrak \FF(L\gfrak)$, 
and acts on $\FF(L\gfrak)_{A\gfrak}$ via the central extension $\hat\theta_R^A$ of $\theta_R^A$ (with the central $c_0$ acting as multiplication by 
$(c+\check{h})^{-1}c\dim\gfrak$).  

These assertions also hold if we replace $\gfrak$ by the abelian Lie algebra $k$ (in particular, $T(\hat D)$ lies in  the closure of $\afrak\circ\lfrakhat$ in $\Ubar\lfrakhat$).
\end{proposition}
\begin{proof}
We only do this for $\gfrak$.
Since $D$ maps $\ann(A)$ to $A\subset L$,
$1\otimes D$ maps the submodule $\gfrak\otimes \ann(A)$ of 
$\gfrak\otimes\omega$ to the submodule $\gfrak\otimes A=A\gfrak$ of 
$\gfrak\otimes L =L\gfrak$. It is clear that  
$\gfrak\otimes \ann(A)$ and $A\gfrak$ are each others 
annihilator relative to the pairing $r_\gfrak$. This implies that
$\hat C (\hat D)$ lies in the closure of the image of 
$A\gfrak\otimes_k \Lg+ \Lg\otimes_k A\gfrak$ 
in $\Ubar \Lghat$. It follows that  $\hat C (\hat D)$ has the form
$cr+\sum_\kappa\sum_{n\ge 1} X_\kappa f_{\kappa,n}\circ X_\kappa g_{\kappa,n}$ with $r\in R$, one of $f_{\kappa_n},  g_{\kappa,n}\in L$  
being in $A$ and the order of $f_{\kappa_n}$ smaller than that of
$g_{\kappa,n}$ for almost all $\kappa, n$. Since the elements of $A$ have 
order $\le 0$ and $X_\kappa f_{\kappa,n}\circ X_\kappa g_{\kappa,n}\equiv
X_\kappa g_{\kappa,n}\circ X_\kappa g_{\kappa,n}\pmod{cR}$, we can assume
that all $f_{\kappa, n}$ lie in $ A$ and so the first assertion follows.

If $D\in \theta_{A,R}$, then for 
$X\in \gfrak$ and $f\in A$, we have $[D,Xf]=X(Df)$, which is an element
of $A\gfrak$ (since $Df\in A$). This shows that $T_\gfrak(\hat D)$ 
preserves $A\gfrak\FF(L\gfrak)$. If $D\in \theta_{A/R}$, then it follows from the proven part  that $T_\gfrak(\hat D)$ maps 
$\FF(L\gfrak)$ to $A\gfrak\FF(L\gfrak)$ and hence induces the zero
map in $\FF(L\gfrak)_{A\gfrak}$. So $\hat\theta_{A,R}$ acts on $\FF(L\gfrak)_{A\gfrak}$ via $\hat\theta_R^A$.
\end{proof}

More relevant than $\FF(L\gfrak)_{A\gfrak}$ will be certain finite dimensional 
quotients thereof  obtained as follows. Let $\ell$ be a positive integer and 
let $V: i\mapsto V_i$ be a map which assigns to every $i\in I$ a 
finite dimensional irreducible representation $V_i$ of $\gfrak$. 
Make $\otimes_{i\in I} V_i$ a $k$-representation of 
$F^0\Lg$ by letting $c$ act as multiplication by 
$\ell$ and $\gfrak\otimes_k\Ocal_i$ on the $i$th factor via its 
projection onto $\gfrak$. If we induce this up to $\Lghat$ we get 
a representation $\widetilde{\HH}_\ell (V)$ of $\Lghat$ 
which clearly is a quotient of $\FF(L\gfrak)$. 
Its irreducible quotient is denoted by $\HH_\ell(V)$. The following is known
(see the book of Kac \cite{kac}). First of all, $\HH_\ell(V)$
is nonzero precisely when each $\gfrak$-representation $V_i$ \emph{is of level $\le\ell$}, i.e., has the property that for every nilpotent  $X\in \gfrak$,  
$X^{\ell +1}$ yields the  zero map in $V_i$. 
Assume this is the case. Then $\HH_\ell(V)$ is integrable as a 
$\Lghat$-module: if $Y\in\gfrak$ is nilpotent and $f\in L$, then $Y f$ acts 
locally nilpotently  in $\HH_\ell(V)$. 

The set of isomorphism classes of finite dimensional 
irreducible $\gfrak$-representations of level $\le\ell$ is finite. 
It is clear that this set, which we denote by $P_\ell$, 
is invariant under dualization (and more generally, under all outer automorphisms of  $\gfrak$). 

Adopting the physicists terminology, we might call 
the $R$-module  of $R$-linear forms on $\HH_\ell(V)_{A\gfrak}$ 
the \emph{conformal block} attached to $A$. The following proposition says 
that it is of finite rank and describes in  essence the WZW-connection.

\begin{proposition}[Finiteness]\label{prop:finiterank}
The space $\HH_\ell (V)$ is finitely generated as a $U\! A\gfrak$-module
(so that $\HH_\ell (V) _{A\gfrak}$ is a finitely generated $R$-module).
The Lie algebra $\hat\theta_R^A$ acts on $\HH_\ell (V) _{A\gfrak}$
via the Sugawara representation with the central 
$c_0$ acting as multiplication by $
\frac{c}{c+\check{h}}\dim\gfrak$. 
\end{proposition}
\begin{proof}
Pick a Cartan subalgebra $\hfrak\subset\gfrak$, let 
$\Delta\subset\hfrak^*$ be
its root system and $\gfrak=\hfrak\oplus_{\alpha\in\Delta}\gfrak_\alpha$
the associated decomposition. Recall that $\gfrak_\alpha$ has a nilpotent
generator and $[\gfrak_\alpha,\gfrak_\beta]$ 
is contained in $\gfrak_{\alpha+\beta}$ when $\alpha+\beta\in\Delta$ and is 
contained in $\hfrak$ otherwise. Recall also that $\gfrak$ is spanned by
$\oplus_{\alpha\in\Delta}\gfrak_\alpha$ and that $\hfrak$ normalizes each
$\gfrak_\alpha$. Choose a generator $t_i$ of $\mfrak_i$. Using a simple 
(PBW-type) argument we see that 
$\HH_\ell (V)$ is spanned over $R$ by the subspaces
\[
(\gfrak_{\alpha_k}q_k)^{\circ r_k}\circ \cdots \circ 
(\gfrak_{\alpha_1}q_1)^{\circ r_1}V
\]
with $k\ge 0$, each $q_\rho$ a negative power of some $t_i$ and such that 
$r_1\ge \cdots \ge r_k$. 
Choose a finite set $Q$ of powers of $t_i$'s that spans an $R$-supplement of  
$\Ocal+A$ in $L$ and let $M$ be the $R$-span of the above subspaces for which all 
$q_\rho$ lie in $Q$. It is clear that then $M$ generates $\HH_\ell (V)$ as a $U\! A\gfrak$-module.
Since each line $\gfrak_\alpha q$ acts locally nilpotently in 
$\HH_\ell(V)$, there exists a positive integer $N$ such that 
$(\gfrak_\alpha q)^{\circ N} V=0$ for all $\alpha\in\Delta$, $q\in Q$. 
So the only generating lines that contribute to $M$ have $r_1< N$ and are hence finite in number. This shows that $M$ is a finitely generated $R$-module. The rest follows from \ref{prop:A}.
\end{proof}

\begin{remark}
We expect the $R$-module $\HH_\ell (V) _{A\gfrak}$ is also flat and that
this is a consequence from a related property for the $UA\gfrak$-module $\HH_\ell (V)$). Such a result, or rather an algebraic proof of it, might simplify the argument 
in \cite{tuy} (see Section \ref{section:doublept}  for our version) which shows that the sheaf of conformal blocks attached 
to a degenerating family of pointed curves is locally free.
\end{remark}

\subsection*{Propagation principle} 
The following proposition is a bare version of what is known as 
\emph{propagation of vacua}, or rather the generalization 
of this fact that can be found in Beauville \cite{beauville}. It often  allows us to reduce the discussion to the case where $I$ is a singleton.

\begin{proposition}\label{prop:singleton}
Let $J\subsetneq I$ be such that  $A$ maps onto $\oplus_{j\in J} L_j/\Ocal_j$. Denote by  $B\subset A$  the kernel, so that we have a surjective Lie homomorphism  $B\gfrak\to\gfrak^{J}$ via which $B\gfrak$ acts on $\otimes_{j\in J} V_j$. Then the map of $B\gfrak$-modules 
$\HH _\ell (V|I-J)\otimes \bigotimes_{j\in J}V_j\to \HH _\ell (V)$
induces an isomorphism  on covariants:
\[\begin{CD}
\left(\HH_\ell (V|I-J)\otimes \bigotimes_{j\in J}V_j\right)_{B\gfrak}
@>\cong>> 
\HH_\ell (V)_{A\gfrak}.
\end{CD}\]
\end{proposition}
\begin{proof} It suffices to do the case when $J$ is a singleton $\{ o\}$. The hypotheses clearly imply  that 
$\HH _\ell (V|I-\{ o\})\otimes V_o\to \HH_\ell (V)_{A\gfrak}$ is onto. 
The kernel is easily shown to be $B\gfrak(\HH _\ell (V|I-\{ o\})\otimes V_o)$. 
\end{proof}

\begin{remark}\label{rem:propvac}
This proposition is sometimes used in the opposite direction:
if $\mfrak_o\subset A$ is a principal ideal with the property that for a generator $t\in\mfrak_o$, the  $\mfrak_o$-adic completion of $A$ gets identified with $R((t))$, then let $\tilde I$ be the disjoint union of $I$ and  $\{ o\}$, $\tilde V$ the extension
of $V$ to  $\tilde I$ which assigns to $o$ the trivial representation and  $\tilde A:=A[t^{-1}]$. With $(\tilde I,\{o\})$ taking the role of $(I,J)$, we  then find that 
$\HH _\ell (V)_{A\gfrak}\cong \HH _\ell (\tilde V)_{\tilde A\gfrak}$.
\end{remark}

\subsection*{Conformal blocks in families}
We transcribe the preceding in more geometric terms. 
This leads us to sheafify many of the notions we
introduced earlier; we shall modify our notation (or its meaning) accordingly. 
Suppose given a proper and  flat morphism between 
$k$-varieties $\pi :\Ccal\to S$ whose base $S$ is smooth and fibers are 
reduced connected curves which have complete intersection
singularities only (we are here not assuming that $\Ccal$ is smooth over $k$).  
Since the family is flat, the arithmetic genus of the fibers is locally constant;
we simply assume it is globally so and  denote this constant genus by $g$.
We also suppose given disjoint sections $S_i\subset\Ccal$, indexed by the 
finite nonempty set $I$  whose union $\cup_{i\in I} S_i$ lies in the smooth part of $\Ccal$ and meets every irreducible component of a 
fiber. The last condition ensures that if $j:\Ccal^\circ:=\Ccal-\cup_{i\in I} 
S_i\subset\Ccal$ is the inclusion, then $\pi j$ is an affine morphism.  

We denote by $(\Ocal_i,\mfrak_i)$ the formal completion of $\Ocal_\Ccal$ along $S_i$, by $\Lcal_i$  the subsheaf of fractions of $\Ocal_i$ with denominator a local generator of  $\mfrak_i$ and by $\Ocal$, $\mfrak$ and $\Lcal$  
the corresponding direct sums. But we keep on using $\omega$, $\theta$, $\hat\theta$ etc.\ for their sheafified counterparts.  So these are now all $\Ocal_S$-modules and the residue pairing is also one of $\Ocal_S$-modules: $r: \Lcal\times \omega\to\Ocal_S$. 
We write $\Acal$ for $\pi_*j_*j^*\Ocal_{\Ccal}$ (a sheaf of $\Ocal_S$-algebras)
and often identify this with its image in $\Lcal$. We denote by $\theta_{\Acal/S}$
the sheaf of $\Ocal_S$-derivations $\Acal\to \Acal$ and by $\omega_{\Acal/S}$
for the direct image on $S$ of the the relative dualizing sheaf of $\Ccal^\circ/S$, 
(in the present situation the relative dualizing sheaf  of $\pi$, $\omega_{\Ccal/S}$, 
is simply the direct image on $\Ccal$ of the sheaf of relative differentials
on the open subset of $\Ccal$ where $\pi$ is smooth). So  $\omega_{\Acal/S}$
is torsion free an hence embeds in $\omega$.

\begin{lemma}
The properties  $A_1$, $A_2$ and $A_3$ hold for the sheaf $\Acal$. Precisely,
\begin{enumerate}
\item[($\Acal_1$)] $\Acal$ is as a sheaf of $\Ocal_S$-algebras flat and of finite type, 
\item[($\Acal_2$)] $\Acal\cap\Ocal=\Ocal_S$ and $R^1\pi_*\Ocal_\Ccal=\Lcal/(\Acal+\Ocal)$ is locally free of rank $g$,
\item[($\Acal_3$)] we have $\theta_{\Acal/S}=\Hom_{\Acal}(\omega_{\Acal/S},\Acal)$ and $\omega_{\Acal/S}$ is the annihilator of $\Acal$ with respect to the residue pairing.
\end{enumerate}
\end{lemma} 
\begin{proof} Property $\Acal_1$ is clear. It is also clear that
$\Ocal_S=\pi_*\Ocal_{\Ccal}\to \Acal\cap \Ocal$ is an isomorphism. 
The long exact sequence defined by the functor $\pi_*$ applied to 
the the short exact sequence
\[
0\to \Ocal_\Ccal\to j_*j^*\Ocal_{\Ccal} \to \Lcal/\Ocal\to 0
\]
tells us that $R^1\pi_*\Ocal_\Ccal=\Lcal/(\Acal+\Ocal)$; in particular,  
the latter is locally free of rank $g$. Hence $\Acal_2$ holds as well.

In order to verify $\Acal_3$, we note that
$\pi_*\omega_{\Ccal/S}$ is the $\Ocal_S$-dual of 
$R^1\pi_*\Ocal_S$, and hence is locally free of rank $g$.
The first part of $\Acal_3$ from the correponding local property
$\theta_{\Ccal/S}=\Hom_{\Ocal_\Ccal}(\omega_{\Ccal/S},\Ocal_\Ccal)$ by apply
$\pi_*j^*$ to either side. This local property is known to hold for families
of curves with complete intersection singularities. (A proof 
under the assumption that $\Ccal$ is smooth---which is does not affect the
generality, since $\pi$ is locally the restriction of that case and both sides
are compatible with base change---runs as follows:
if $j':\Ccal'\subset \Ccal$ denotes the locus where $\pi$ is smooth, then
its complement is of codimension $\ge 2$ everywhere.
Clearly, $\theta_{\Ccal/S}$ is the $\Ocal_{\Ccal}$-dual of 
$\omega_{\Ccal/S}$ on $\Ccal'$ and since both are inert under $j'_*j'{}^*$, 
they are equal everywhere.) 

The last assertion  essentially restates the well-known fact that
the polar part of a rational section of $\omega_{\Ccal/S}$
must have zero residue sum, but can otherwise be arbitrary. More precisely, 
the image of $\omega_{\Acal/S}$ in $\omega/F^1\omega$ is the 
kernel of the residue map $\omega/F^1\omega\to \Ocal_S$. 
The intersection $\omega_{\Acal/S}\cap F^1\omega$ is $\pi_*\omega_{\Ccal/S}$ 
and is hence  locally free of rank $g$.
Since $\ann (F^1\omega)=\Ocal$, it follows that 
$\ann(\omega_{\Acal/S})\cap \Ocal$ and $\Lcal/(\ann(\omega_{\Acal/S})+\Ocal)$ 
are locally free of rank 1 and  $g$ respectively. Since $\Acal$ has these properties also 
and is contained in $\ann(\omega_{\Acal/S})$, we must have 
$\Acal=\ann(\omega_{\Acal/S})$.
\end{proof}

For what follows one usually supposes 
that the fibers are stable $I$-pointed curves (meaning that  every fiber of $\pi j$
has  only ordinary double points as singularities and has finite 
automorphism group) and is versal (so that in the discriminant $D_\pi$ of 
$\pi$ is a reduced normal crossing divisor), but we shall not make these 
assumptions yet. Instead, we assume the considerable weaker property that 
the sections of the sheaf $\theta_S(\log D_\pi)$ of vector fields on $S$ 
tangent to $D_\pi$ lift locally on $S$ to vector fields on $\Ccal$. 
This is for instance the case  if $\Ccal$ is smooth and $\pi$ is multi-transversal with respect 
to the (Thom) stratification of $\Hom (T\Ccal,\pi^*TS)$ by rank
\cite{looij:book}. Notice that the restriction homomorphism
$\theta_S(\log D_\pi)\otimes \Ocal_{D_\pi}\to \theta_{D_\pi}$ is
an isomorphism.

Let $\theta_{\Ccal, S}\subset\theta_\Ccal$ 
denote the sheaf of derivations which preserve $\pi^*\Ocal_S$.
Applying $\pi_*j_*j^*$ to the exact sequence 
$0\to\theta_{\Ccal/S}\to\theta_{\Ccal, S}\to 
\theta_{\Ccal, S}/\theta_{\Ccal/S}\to 0$ yields  the exact sequence
\[
0\to \theta_{\Acal}\to\theta_{\Acal, S}\to \theta_S(\log D_\pi)\to 0.
\] 
(using our litfability assumption and the fact that $\pi j$ is affine).

We have defined $\hat\theta_{\Acal,S}$ as the preimage of  
$\theta_{\Acal,S}$ in $\hat\theta_{\Lcal,S}$ and 
$\hat\theta_S(\log D_\pi)$ as the quotient 
$\hat\theta_{\Lcal ,S}/\theta_\Acal$. 
These centrally extend $\theta_{\Acal,S}$  and  
$\theta_S$ by $c_0 \Ocal_S$.

Observe that $\Lcal\gfrak=\gfrak\otimes_{\Ocal_S}\Lcal$ is now a sheaf of 
Lie algebras over  $\Ocal_S$. The same applies to $\hat\lfrak$ and so we have a Virasoro 
extension $\hat\theta_S$ of $\theta_S$ by $c_0\Ocal_S$. We have also defined
${\Acal\gfrak}= \gfrak\otimes_k \Acal$, which is  a Lie subsheaf of 
$\Lcal\gfrak$ as well as of $\Lcalhatg$ and the Fock type $\Lcalhatg$-module
$\Fcal(\Lcalhatg)$. The \emph{universal covacuum module} is the 
sheaf of ${\Acal\gfrak}$-covariants in the latter, 
\[
\Fcal(\Lcalhatg)_{\Ccal}:=\Fcal(\Lcalhatg)_{\Acal\gfrak}:=\Fcal(\Lcalhatg)/
{\Acal\gfrak}\Fcal(\Lcalhatg).
\]
From Proposition \ref{prop:A} we get:

\begin{corollary}\label{cor:thetaaction}
The representation of the Lie algebra $\hat\theta_{\Acal,S}$ on 
$\Fcal(\Lcalhatg)$ preserves  ${\Acal\gfrak}\Fcal(\Lcalhatg) $
and acts on $\Fcal(\Lcalhatg) _{\Ccal}$ via the central extension 
$\hat\theta_S(\log D_\pi)$ of $\theta_S(\log D_\pi)$ 
(with the central $c_0$ acting as multiplication by 
$(c+\check{h})^{-1}c\dim\gfrak$). This construction has a base change property
along the smooth part $D_\pi^\circ$ of the discriminant: we may identify
$\Fcal(\Lcalhatg) _{\Ccal}\otimes\Ocal_{D_\pi^\circ}$ with 
$\Fcal(\widehat{\Lcal|_{D_\pi^\circ}\gfrak})_{\Ccal |_{D_\pi^\circ}}$ and 
the action of $\hat\theta_S(\log D_\pi)\otimes\Ocal_{D_\pi^\circ}$ on
it factors through $\hat\theta_{D_\pi^\circ}$.   
\end{corollary}

The bundle of integrable representations $\Hcal_\ell (V)$ over $S$ is 
defined in the expected manner: it is obtained as a quotient of 
$\Fcal(\Lcalhatg)$ in the way $\HH_\ell (V)$ is obtained from $\FF(L\gfrak)$.
According to Corollary \ref{cor:thetaaction}, $\Hcal_\ell(V)$ comes with a 
$k$-linear representation of the Lie algebra $\hat\theta_S(\log D_\pi)$ that 
satisfies the Leibniz rule and on which $c_0$ acts as 
multiplication by $\ell (\ell +\check{h})^{-1}\dim\gfrak$.
We write $\Hcal_\ell(V)_{\Ccal}$ for $\Hcal_\ell(V)_{\Acal\gfrak}$.

\begin{corollary}[WZW-connection]\label{cor:flatter}
The $\Ocal_S$-module $\Hcal_\ell(V)_{\Ccal}$ is of finite rank;  
It is also locally free over $S-D_\pi$ and the Lie action of 
$\hat\theta_S(\log D_\pi)$  on  $\Hcal_\ell (V)_{\Ccal}$ 
defines a flat connection on the associated bundle of 
projective spaces, $\PP_S(\Hcal_\ell (V)_{\Ccal})\to S$ with a 
logarithmic singularity along $D_\pi$. The same base 
change property holds along the smooth part $D_\pi^\circ$ of the 
discriminant as in Corollary \ref{cor:thetaaction}.
\end{corollary}  
\begin{proof}
The first assertion follows from \ref{prop:finiterank}. 
The rest is clear except perhaps the 
local freeness of $\Hcal_\ell(V)_{\Ccal}$ on $S-D_\pi$. But this follows 
from the local
existence of a connection in the $\Ocal_S$-module $\Hcal_\ell(V)_{{\Ccal}}$.      
\end{proof}

\begin{remark}\label{rem:topmeaning}
In case the base field $k$ is $\CC$ we could equally well work in the complex-analytic  category. The preceding corollary then leads to an interesting topological functor:
if $\Sigma$ is a compact smooth oriented surface, $I\subset\Sigma$ a finite nonempty subset such that for every $i\in I$ is given a $\gfrak$-representation $V_i$, then there is canonically associated to these data a projective space  $\PP(\Sigma,I,V)$: for every complex structure on $\Sigma$ compatible with the
given smooth structure and orientation we have defined a conformal block as above.
The complex structures with these property form a simply connected (even contractible) space $\Ccal (\Sigma)$. To see this:  such a structure amounts to giving a conformal structure on $\Sigma$, that is, a Riemann metric up to multiplication by the exponential of a real valued function. The space of Riemann metrics is a convex cone in some linear space (hence contractible) and the linear space of real valued functions on $\Sigma$ acts sufficiently nice on it to ensure that the orbit space remains simply connected. So any two such complex structures can be joined a path whose homotopy class is unique. Corollary \ref{cor:flatter} then implies that the WZW-connection enables us to identify their associated conformal blocks in a unique manner. Perhaps the best way to describe $\PP(\Sigma,I,V)$ is as the space of flat sections of the flat projective space
bundle of conformal blocks over $\Ccal (\Sigma)$. The naturality implies that the group of orientation preserving diffeomorphisms  of the pair $\Sigma$ which fix $I$
pointwise acts on $\PP(\Sigma,I,V)$. Since its identity component is contractible, the
action will be through the mapping class group $\G(\Sigma,I)$. It is conjectured that this action leaves invariant a Fubini-Study metric on $\PP(\Sigma,I,V)$. 
\end{remark}

\subsection*{Propagation principle continued} In the preceding subsection we made the 
assumption throughout that a union $\tilde S$ of sections of $\Ccal\to S$ is given
to ensure that $\Ccal-\tilde S\to S$ is affine. However, the propagation principle
permits us to abandon that assumption. In fact, this leads us to let $V$ stand for any map which  assigns to every closed point $p$ of $\Ccal$ an irreducible $\gfrak$-representation 
$V_p$ of level $\le \ell$, subject to the condition  that  its \emph{support}, $\supp (V)$ (i.e., the set of $p$ for which $V_p$ 
is not the trivial representation),  is a trivial finite cover over $S$ and contained in the locus where $\pi: \Ccal\to S$ is smooth.
Our earlier defined $\Hcal_\ell(V)$ is now better written as 
$\Hcal_\ell(V|_{\supp (V)})$. Since  $\Ccal-\supp (V)$ need not be 
affine over $S$, this does not yield the right notion of conformal block. 
We can find however, at least locally over $S$, additional sections $S'_i$
of  $\Ccal-\supp (V)\to S$ so that $\Ccal^\circ:=\Ccal-\supp (V)-\cup_i S'_i$
is affine over $S$. Then we can form $\Hcal_\ell(V|\Ccal-\Ccal^\circ)$ and 
Remark  \ref{rem:propvac}  shows that the resulting conformal block 
$\Hcal_\ell(V|\Ccal-\Ccal^\circ)_{\gfrak[\pi_*\Ocal_{\Ccal^\circ}]}$ is 
independent of the choices made.  This suggests that we let $\Hcal_\ell(V)$ 
stand for the sheaf associated to the presheaf
\[
S\supset U\mapsto \varprojlim_{\tilde S} \Hcal_\ell(V|_{\tilde S}),
\]
where $\tilde S$ runs over the closed subvarieties of $\pi^{-1}U$ that trivial finite 
covers of $S$, contain $\supp (V)$ and have affine  complement in $\pi^{-1}U$ 
and perhaps also justifies our custom of writing $\Hcal_\ell (V)_\Ccal$ 
for the associated sheaf of conformal blocks.  
It is clear that in this set-up there is also no need anymore to insist 
that the fibers of $\pi$ be connected.

\subsection*{The genus zero case and the KZ-connection}
We take $C\cong\PP^1$ and let $\ell$ and $z\in C\mapsto V_z$  be as usual.
Choose an affine coordinate $z$ on $C$  such that $\infty\notin\supp(V)$ and let $z_1,\dots ,z_N\in\CC$ enumerate  the distinct points of $\supp (V)$. Write $V_i$ for $V_{z_i}$. The local field  attached to $\infty$ has parameter $t_\infty=z^{-1}$. If $\HH_\ell$ denotes the representation of $\widehat{\gfrak((z^{-1}))}$ attached to the trivial  representation (with generator $v_0$), then by the propagation principle 
\ref{prop:singleton} we have $\HH_\ell (V)_C=(\HH_\ell\otimes V_1\otimes
\cdots\otimes V_N)_{\gfrak[[z]]}$, where $\gfrak[[z]]$ acts on $V_i$ via its 
evaluation at $z_i$. According to \cite{kac}, the $\gfrak[[z]]$-homomorphism 
$U\gfrak[[z]]\to \HH_\ell$, $u\mapsto uv_0$ is surjective and its kernel is the 
left ideal generated by the constants 
$\gfrak$ and $(Xz)^{1+\ell}$, where $X\in\gfrak$ is a regular nilpotent element
(as these form a single orbit under the adjoint representation, it does not matter
which one we take).  This implies that $\HH_\ell (V)_{\PP^1}$ can be identified with 
a quotient of the space of $\gfrak$-covariants 
$(V_1\otimes\cdots\otimes V_N)_{\gfrak}$, namely its biggest quotient 
on which  $(\sum_{i=1}^N z_iX_{(i)})^{1+\ell}$ acts trivially (where $X_{(i)}$ acts on
$V_i$ as $X$ and on the other tensor factors $V_j$, $j\not= i$, as the identity).

Now regard $z_1,\dots ,z_N$ as variables. If $S$ denotes the 
open subset of $(z_1,\dots ,z_N)\in\CC^N$ of pairwise distinct $N$-tuples, then
we have in an evident manner a family $\Ccal$ over $S$,  $\pi: C\times S\to S$, with $N+1$ sections (including the one at infinity) so that we also have defined
$\Ccal^\circ\to S$. This determines a sheaf of conformal blocks $\Hcal_\ell (V)_\Ccal$ over $S$. According to the preceding, we have an exact sequence
\[
\begin{CD}
(V_1\otimes\cdots\otimes V_N)_{\gfrak}\otimes_k\Ocal_S @>{(\sum_i z_iX_{(i)})^{1+\ell}}>>
(V_1\otimes\cdots\otimes V_N)_{\gfrak}\otimes_k\Ocal_S \to \Hcal_\ell (V)_\Ccal\to 0.
\end{CD}
\]
We  identify  its WZW connection, or rather, a natural lift of that connection to
$V_1\otimes\cdots\otimes V_N\otimes_k\Ocal_S$. 
In order to compute the covariant derivative with respect to the vector field 
$\p_i:=\frac{\p}{\p z_i}$ on $S$, we follow our receipe and lift it
to $\Ccal$ in the obvious way (with zero component along $C$). We continue to denote that lift by $\p_i$ and determine its (Sugawara) action on  $\Hcal_\ell (V)$.
We first observe that $\p_i$ is tangent to all the sections, except the $i$th. Near that section we decompose it as $(\frac{\p}{\p z}+\p_i)-\frac{\p}{\p z}$, where the first term
is tangent to the $i$-th section and the second term is vertical. The action of the former is easily understood (in terms of an obvious trivialization of $\Hcal_\ell (V)$, it acts as derivation with respect to $z_i$). The vertical term, $-\frac{\p}{\p z}$, acts via
the Sugawara representation, that is, it acts on the $i$th slot as 
$-\frac{1}{\ell +\check{h}}\sum_\kappa X_\kappa(z-z_i)^{-1}\circ X_\kappa$ and  as the identity on the others. This action does not induce one in 
$V_1\otimes\cdots\otimes V_N\otimes_k\Ocal_S$. To make it so, we add to this the action of the element 
\[
\left(\frac{1}{\ell +\check{h}}\sum_\kappa X_\kappa(z-z_i)^{-1}\right)\circ X^{(i)}_\kappa 
\in  \gfrak[\Ccal^\circ]U\Lcalhatg, 
\]
for this sum then acts in $V_1\otimes\cdots\otimes V_N\otimes_k\Ocal_S$ as 
\[
\frac{1}{\ell+\check{h}}\sum_{j\not=i} \frac{1}{z_j-z_i}X_\kappa^{(j)}X_\kappa^{(i)}.
\]
Let us regard the Casimir element $c$ as a tensor of $\gfrak\otimes_k\gfrak$, 
and denote by $c^{(i,j)}$ its action in $V_1\otimes\cdots\otimes V_N$ on the $i$th and $j$th factor (since $c$ is symmetric, we have $c^{(i,j)}=c^{(j,i)}$, so  that we need not worry about the order here).
We conclude that the WZW-connection is induced by the connection
on $V_1\otimes\cdots\otimes V_N\otimes_k\Ocal_S$ defined by the connection form
\[
\frac{1}{\ell+\check{h}}\sum_{i=1}^N \sum_{j\not=i} \frac{dz_i}{z_j-z_i}c^{(i,j)}=
-\frac{1}{\ell+\check{h}}\sum_{1\le i<j\le N} \frac{d(z_i-z_j)}{z_i-z_j}c^{(i,j)}
\]
This lift of the WZW-connection  is known as the Knizhnik-Zamolodchikov connection. It is not difficult to verify that it is still flat (see for instance \cite{kohno}).

\section{Double point degeneration}\label{section:doublept}

\subsection*{Factorization} 
Let $\ell$, $\pi :\Ccal\to S$ and $V$ be as in the preceding section so that
we have defined the sheaf of conformal blocks $\Hcal_\ell (V)_\Ccal$.
We assume in addition that we are given a section $S_0$ along which $\pi$ 
has an ordinary double point and a partial  normalization 
$\nu: \tilde\Ccal\to\Ccal$ which only separates the branches. 
So $\nu$ is an isomorphism over the complement of $S_0$ 
and $S_0$ has two disjoint lifts  to $\Ccal$, which we shall denote by $S_+$ and $S_-$. The \emph{factorization principle} expresses $\Hcal_\ell (V)_\Ccal$ in  terms of conformal blocks attached to the partial normalization of $\tilde\Ccal$.

Recall that $P_\ell$ denotes 
the set of isomorphism classes of irreducible $\gfrak$-representations 
of level $\le \ell$ and is invariant under dualization: if 
$\lambda\in P_\ell$, then $\lambda^*\in P_\ell$. 
Let  $V_\lambda$  be a $\gfrak$-representation in the equivalence class 
$\lambda\in P_\ell$ and choose  $\gfrak$-equivariant dualities 
$b_\lambda: V_\lambda\otimes V_{\lambda^*}\to k$, $\lambda\in P_\ell$.
 
\begin{proposition}\label{prop:fact}
Let $\tilde V_{\lambda,\lambda^*}$  be the representation valued map on $\tilde\Ccal$ which is constant equal to $V_\lambda$ resp.\
$V_{\lambda^*}$ on $S_+$  resp.\ $S_-$ and on any other closed point of
the representation already assigned to its image in $\Ccal$.
Then the contractions $b_\lambda$ define an isomorphism
\[
\begin{CD}
\oplus_{\lambda\in P_\ell} \Hcal _\ell (\tilde V_{\lambda,\lambda^*})_{\tilde \Ccal}
@>\cong>>\Hcal _\ell (V)_\Ccal.
\end{CD}
\]
\end{proposition}

This is almost a formal consequence of:

\begin{lemma}\label{lemma:bilevel}
Let $M$ be a finite dimensional representation of $\gfrak\times\gfrak$
which is  of level $\le \ell$ relative to both factors. If $M^\delta$ denotes 
that same space viewed as $\gfrak$-module with respect to the 
diagonal embedding $\delta: \gfrak\to\gfrak\times\gfrak$, then the 
contraction
$\oplus _{\lambda\in P_\ell} M\otimes (V_\lambda\boxtimes V_\lambda^*)\to M$ 
(each component is defined  by $b^\lambda$;  the symbol $\boxtimes$ 
stands for the exterior tensor product of representations) 
induces an isomorphism between  covariants: 
\[
\begin{CD}
\oplus _{\lambda\in P_\ell} \left(M\otimes (V_\lambda\boxtimes 
V_\lambda^*)\right)_{\gfrak\times\gfrak}@>\cong>> M^\delta_\gfrak.
\end{CD}
\]
\end{lemma}
\begin{proof}
Without loss of generality we may assume that $M$ is irreducible, or more precisely, 
equal to  $V_\mu\boxtimes V_{\mu'}$ for some $\mu,\mu'\in P_\ell$.
Then $M^\delta=V_\mu\otimes V_{\mu'}$. By Schur's lemma, $M^\delta_\gfrak$ is one-dimensional  if $\mu'=\mu^*$ and trivial otherwise. That same lemma applied
to $\gfrak\times\gfrak$ shows that 
$(M\otimes (V_\lambda\boxtimes V^*_\lambda))_{\gfrak\times\gfrak}$  
is zero unless $(\mu,\mu')=(\lambda^*,\lambda)$, in which case it is one-dimensional. The lemma follows.
\end{proof}

\begin{proof}[Proof of \ref{prop:fact}] The issue is local on $S$ and so we may assume
that there exists an  an affine open-dense subset $\Ccal^\circ$ of
$\Ccal$ which contains its singular locus, is disjoint with $\supp (V)$ and
is such that its complement $\tilde S:=\Ccal-\Ccal^\circ$ is a trivial finite cover 
over $S$.
Then $\Hcal _\ell (V)_\Ccal=\Hcal_\ell (V|\tilde S)_{\pi_*\Ocal_{\Ccal^\circ}[\gfrak]}$.  Let $\tilde\Ccal^\circ:=\nu^{-1}\Ccal^\circ$.  
Evaluation in $S_0$ resp.\ 
$S_+,S_-$ defines epimorphisms $\pi_*\Ocal_{\Ccal^\circ}\to \Ocal_S$ resp.\ $\pi_*\nu_*\Ocal_{\tilde\Ccal^\circ}\to \Ocal_S\oplus  \Ocal_S$ whose kernels may be identified by means  of $\nu$. If we denote that common kernel by $\Ical$ and $\Ical\gfrak$ has the evident 
meaning, then the argument used to prove Proposition \ref{prop:finiterank} 
shows that $M:=\Hcal_\ell (V|_{\tilde S})_{\Ical\gfrak}$ is a $\Ocal_S$-module of finite rank. It underlies a representation of $\Ocal_S\gfrak\oplus\Ocal_S\gfrak$ of level $\le \ell$ relative to both factors and is such that $M^\delta_\gfrak=\Hcal_{\ell}(V)_\Ccal$. Now apply
Lemma \ref{lemma:bilevel}.
\end{proof}

\subsection*{Local freeness}
We continue with the situation of Proposition \ref{prop:fact}, except that we now 
assume (for simplicity only, actually) that the base $S$ is $\spec (k)$: $C$ is a reduced complete curve with 
complete intersection singularities only which has at $S_0\in C$ an ordinary double
point. Choose generators 
$t_{\pm}$ of the maximal ideals of the completed local rings of $C$ at  $S_{\pm}$. 
There is then a canonical \emph{smoothing} of $C$, that is, a way of making $C$ the closed fiber of a flat morphism $\tau: \Ccal\to \Delta$, with $\Delta$  the spectrum of the discrete valuation ring $k[[\tau]]$, such that  
the generic fiber is smooth: in the product $\tilde C\times\Delta$, blow up 
$(S_{\pm}, o)$ and let $\tilde\Ccal$ be the formal neighborhood of the strict transform of $\tilde C\times\{ o\}$. 
So at the preimage of  $(S_{\pm}, o)$ on the strict transform of $\tilde C\times\{ o\}$ we have the formal 
coordinate chart $(t_{\pm}, \tau/t_{\pm})$. Now let $\Ccal$ be the quotient of $\tilde\Ccal$ obtained by identifying these formal charts up to order: $(t_+, \tau/t_+)=(\tau/t_{-},t_-)$, so that $(t_+,t_-)$ is now a formal chart  of $\Ccal$ on which we have 
$\tau=t_+t_-$ (in either domain  $\tau$ represents the same regular function).  We thus have defined a flat morphism 
$\tau: \Ccal\to \Delta$ whose closed fiber may be identified with $C$. The domain $\Ccal$ is smooth over $k$ and that the generic fiber of $\tau$ is smooth over $k((\tau ))$. This is our canonical smoothing of $C$. 
Let us also notice that $\Ccal$ is at every $p\in C-\{S_0\}$  canonically identified with 
$(C,p)\times\Delta$ with $\tau$ given as the projection on the second factor.  

Let $I$ be a finite subset of the smooth part of $C$ such that $C-I$ is affine and 
$V$ is a $\gfrak$-representation valued map on $I$. We extend $V$ 
canonically to $\Ccal$  by letting it be zero at $S_0$ and constant 
on $\{p\}\times\Delta$ when $p\not= S_0$; we denote that extension 
by $V\!/\!\Delta$.
Then we have defined the $k[[\tau]]$-module $\Hcal_\ell (V\!/\!\Delta)$. 
It is clear that $\Hcal_\ell (V\!/\!\Delta)$ is naturally identified with 
$\HH_\ell (V)[[\tau]]$. We assume that the complement 
of the support of $V$ is affine and we denote by $A$ its $k$-algebra of regular
functions. Then the complement of the support of 
$V\!/\!\Delta$ is affine over $\Delta$ and if $\Acal$ denotes the corresponding 
$k[[\tau]]$-algebra of regular functions, then $A=\Acal/\tau\Acal$.
According to Proposition \ref{prop:finiterank}, the conformal block 
$\Hcal_\ell (V\!/\!\Delta)_{\Ccal}$ is a finitely generated $k[[\tau]]$-module; its reduction 
modulo $\tau$ is clearly $\HH_\ell (V)_{A\gfrak}$. If  we denote by
$B$ the algebra of regular functions on $C-\supp(V)-S_0=\tilde C-
\supp(V)-(S_+ \cup S_-)$, then Proposition 
\ref{prop:fact} identifies the latter with 
$\oplus_{\lambda\in P_\ell} \HH _\ell (\tilde V_{\lambda,\lambda^*})_{B\gfrak}$.
It is our goal to extend this identification to one of 
$\Hcal_\ell (V\!/\!\Delta)_{\Ccal}$ with
$\oplus_{\lambda\in P_\ell} \HH _\ell (\tilde V_{\lambda,\lambda^*})_{B\gfrak} 
[[\tau]]$ (which will imply that $\Hcal_\ell (V\!/\!\Delta)_{\Ccal}$ is a free
$k[[\tau]]$-module). 

Put $\Ocal_*:=k[[t_+]]\oplus k[[t_-]]$ and $L_*:=k((t_+))\oplus k((t_-))$.  

\begin{lemma}\label{lemma:expansion}
The rational map $\tilde C\times\Delta\leftarrow \tilde\Ccal\to \Ccal$ identifies 
$k[[t_+,t_-]]$ with the subalgebra of  $L_*[[\tau]]$ of elements of the form 
$\sum_{n\ge 0,m\ge 0}a_{m,n}(t_+^{m-n}\tau^n,t_-^{n-m}\tau^m)$.
Furthermore, any  continuous $k$-derivation $k[[t_+,t_-]]\to k[[t_+,t_-]]$ which fixes
$\tau=t_+t_-$ defines in $L_*$ an operator of the form
\[
(t_+\frac{\p}{\p t_+}, 0) + \sum_{m,n\ge 0} a_{m,n}(t_+^{m-n+1}\tau^n\frac{\p}{\p t_+},
-t_-^{n-m+1}\tau^m\frac{\p}{\p t_-}).
\]  
\end{lemma}
\begin{proof}
Let $f=\sum_{n,m\ge 0} a_{m,n}t_+^mt_-^n$.
If we substitute $t_-=\tau/t_+$, this becomes at $(S_+,o)$: 
$\sum_{n\ge 0} (\sum_{m\ge 0}a_{m,n}t_+^{m-n})\tau^n$. So
the coefficient of $\tau^n$ has polar part $\sum_{m=0}^{n-1}a_{m,n}t_+^{m-n}$ and constant term $a_{n,n}$. Similarly at
$(S_-,o)$, $f$ is written as $\sum_{m\ge 0} (\sum_{n\ge 0}a_{m,n}t_+^{n-m})\tau^m$
with the coefficient of $\tau^m$ having polar part $\sum_{n=0}^{m-1}a_{m,n}t_+^{n-m}$ and constant term $a_{m,m}$. The polar parts and constant terms of these expressions determine all the $a_{n,m}$.

The second assertion is left as an exercise. 
\end{proof}

Given $\lambda\in P_\ell$, then the Casimir element $c=\sum_\kappa X_\kappa\circ X_\kappa$ acts in 
$V_\lambda$ as a scalar, a scalar we shall denote by 
$c(\lambda)$. Observe that $c(\lambda^*)=c(\lambda)$.

Let $\HH^{\pm}_\ell(V_\lambda)$ denote the 
representation attached to $V_\lambda$ of the central extension
of $\gfrak((t_\pm))$, so that  
$\HH^+_\ell(V_\lambda)\otimes \HH^-_\ell(V_{\lambda^*})$
is a representation of the central extension $\widehat{L_*\gfrak}$  of $L\gfrak$.

\begin{lemma}\label{lemma:invariant}
There exists a tensor valued Laurent series
\[
\varepsilon^\lambda=\sum_{d=0}^\infty \varepsilon^\lambda_d\tau^d\in (\HH^+_\ell(V_\lambda)\otimes 
\HH^-_\ell(V_{\lambda^*}))[[\tau]]\; \text{  with $\varepsilon^\lambda_d\in
F^{-d}\HH^+_\ell(V_\lambda)\otimes F^{-d}\HH^-_\ell(V_{\lambda^*})$}
\]
whose constant term $\varepsilon^\lambda_0\in V_\lambda\otimes V_{\lambda^*}$ is the dual of $b_\lambda$ and which 
is annihilated by the image of $\gfrak[[t_+,t_-]]$ in $\widehat{L_*\gfrak}$. Moreover, for any continuous $k$-derivation 
$D: k[[t_+,t_-]]\to k[[t_+,t_-]]$ which fixes $\tau$,  $\varepsilon^\lambda$ is an eigenvector of $T_\gfrak(D)$
with eigenvalue $-\frac{c(\lambda)}{2(\ell+\check{h})}$.
\end{lemma}
\begin{proof}
We first observe that the choice of the coordinates $t_+$ and $t_-$ defines 
a grading on all the relevant objects on which we have defined the associated
filtration $F$
(e.g., the degree zero summand of   $\HH_\ell (V_\lambda)$ is $V_\lambda$). It is known \cite{kac} that the  pairing $b_\lambda: V_\lambda\times V_{\lambda^*}\to k$ extends (in a unique manner) to a perfect pairing  
\[
b_\lambda : \HH^+_\ell(V_\lambda)\times
\HH^-_\ell(V_{\lambda^*})\to k
\]
characterized by the property that 
$b_\lambda(Xt_+^k u,u')+b_\lambda(u,Xt_-^{-k} u')=0$ for all $X\in\gfrak$ 
and $k\in\ZZ$. It follows from that the restriction of
$b_\lambda$ to  
$\HH^+_\ell(V_\lambda)_{-d}\times \HH^-_\ell(V_{\lambda^*})_{-d'}$ is zero 
when $d\not= d'$ and is perfect when $d=d'$.  So if $\varepsilon^\lambda_d\in
\HH^+_\ell(V_\lambda)_{-d}\otimes \HH^-_\ell(V_{\lambda^*})_{-d}$ denotes the latter's inverse, then we have for all $k\in\ZZ$, $X\in\gfrak$
\[
(Xt_+^k,0) \varepsilon^\lambda_{d+k} +(0,Xt_-^{-k})\varepsilon^\lambda_{d}=0.
\]
This amounts to the property that 
$(Xt_+^k, \tau^kXt_-^{-k})$ annihilates  $\sum_{d\ge 0}\varepsilon_d^\lambda \tau^d$. Since Lemma \ref{lemma:expansion} says that the elements of $\gfrak 
[[t_+,t_-]]$ have series expansions in  $(Xt_+^k, Xt_-^{-k})\in 
\widehat{L_*\gfrak}$ with coefficients in $k[[\tau]]$, the first statement of the lemma follows.

The second will be a straightforward computation. Using Lemma \ref{lemma:expansion},
we write  $D$ as an operator in $L_*$.  We then find that we need to verify the following two assertions:
 \begin{enumerate}
\item[(i)] For every pair $m,n\ge 0$, $\tau^nT_\gfrak(t_+^{m-n+1}\frac{\p}{\p t_+})
-\tau^mT_\gfrak(t_-^{n-m+1}\frac{\p}{\p t_-})$ kills $\varepsilon^\lambda$ and
\item[(ii)] $T_\gfrak (t_+\frac{\p}{\p t_+})(\varepsilon^\lambda)=
-\frac{c(\lambda)}{2(\ell+\check{h})}\varepsilon^\lambda$.
 \end{enumerate} 
As to (i), if we substitute 
\[
T_\gfrak(t_+^{m-n+1}\frac{\p}{\p t_+})=-\frac{1}{2(\ell+\check{h})}
\sum_{j\in\ZZ}\sum_\kappa :X_\kappa t_+^{m-n-j}\circ X_\kappa t_+^j:
\]
 and do likewise
for $T_\gfrak(t_-^{n-m+1}\frac{\p}{\p t_-})$, this assertion follows easily. 

For (ii) we first observe that $T_\gfrak (t_+\frac{\p}{\p t_+})$ 
preserves the grading of
$\HH^+_\ell(V_\lambda)$ and acts on $\HH^+_\ell(V_\lambda)_0=V_{\lambda}$ as 
$-(2\ell+2\check{h})^{-1}\sum_\kappa X_\kappa\circ X_\kappa$. This is just 
multiplication by $-\frac{c(\lambda)}{2(\ell+\check{h})}$. 
If $u\in \HH^+_\ell(V_\lambda)_{-d}$ is written 
$u=Y_rt_+^{-k_r}\circ\cdots\circ Y_1t_+^{-k_1}\circ v$ with $v\in V_\lambda$, $Y_\rho\in \gfrak$, $d=k_r+\cdots +k_1$, then 
$T_\gfrak (t_+\frac{\p}{\p t_+})(u)=-du+Y_rt_+^{-k_r}\circ\cdots\circ Y_1t_+^{-k_1}\circ
T_\gfrak (t_+\frac{\p}{\p t_+})(v)=(-d-\frac{c(\lambda)}{2(\ell+\check{h})})u$.
Since $t_+\frac{\p }{\p t_+}\tau^d=d\tau^d$, it follows that $\varepsilon^\lambda_d\tau^d$ is an eigenvector
of  $T_\gfrak (t_+\frac{\p}{\p t_+})$ with eigenvalue $-\frac{c(\lambda)}{2(\ell+\check{h})}$.
\end{proof}

\begin{theorem}\label{thm:doublept}
The $k[[\tau]]$-homomorphism  
\begin{align*}
E=(E_\lambda)_\lambda: \Hcal_\ell (V\!/\!\Delta)=\HH_\ell (V)[[\tau]]&\to
\oplus_{\lambda\in P_l} \HH_\ell (\tilde V_{\lambda,\lambda^*})[[\tau]],\\
U=\sum_{k\ge 0} u_k\tau^k&\mapsto \sum_{\lambda\in P_\ell} U\varepsilon^\lambda :=
\sum_{\lambda\in P_\ell}\sum_{k,d\ge 0} u_k\otimes \varepsilon_d^\lambda\tau^{k+d}
\end{align*}
is also a map of $\Acal\gfrak$-representations if we let $\Acal\gfrak$ act on the
left hand side via the inclusion
$\Acal\subset B[[\tau]]$. The resulting  $k[[\tau]]$-homomorphism of covariants,
\[
\Hcal_\ell (V\!/\!\Delta)_{\Ccal}\cong 
\oplus_{\lambda\in P_l} \HH_\ell (\tilde V_{\lambda,\lambda^*})_{B\gfrak}[[\tau]]
=\oplus_{\lambda\in P_l} \HH_\ell (\tilde V_{\lambda,\lambda^*})_{\tilde C}[[\tau]]
\]
is an isomorphism. In particular, $\Hcal_\ell (V\!/\!\Delta)_{\Ccal}$ is a free 
$k[[\tau]]$-module. Moreover, covariant derivation in 
$\Hcal_\ell (V\!/\!\Delta)_{\Ccal}$ by
$\tau\frac{d}{d\tau}$ respects each summand  
$\HH_\ell (\tilde V_{\lambda,\lambda^*})_{\tilde C}[[\tau]]$ 
and acts there as the first order differential operator 
$\tau\frac{d}{d\tau}+\frac{c(\lambda)}{2(\ell+\check{h})}$.
\end{theorem}
\begin{proof}
The first statement is immediate from Lemma \ref{lemma:invariant}.
So the map on covariants is defined and is  $k[[\tau]]$-linear. If we reduce modulo $\tau$ we  get the map 
\[
u\in \HH_\ell (V)_{A\gfrak}\mapsto \sum_{\lambda\in P_\ell} u\otimes \varepsilon_0^\lambda\in 
\oplus_{\lambda\in P_l} \HH_\ell (\tilde V_{\lambda,\lambda^*})_{B\gfrak}.
\]
We recognize its domain and range as $\Hcal_\ell (V\!/\!\Delta)_C$ and  
$\oplus_{\lambda\in P_l} \HH_\ell (\tilde V_{\lambda,\lambda^*})_{\tilde C}$
and we observe that the map itself is just the inverse of the isomorphism of Proposition \ref{prop:fact}. Since the range is a free $k[[\tau]]$-module, this implies that the 
$k[[\tau]]$-homomorphism of covariants is an isomorphism.

According to Corollary \ref{cor:derive} covariant derivation with respect to $\tau\frac{d}{d\tau}$ in $\Hcal_\ell (V\!/\!\Delta)_{\Ccal}$
is defined by means of a $k$-derivation $D$ of  $\Acal$ which lifts $\tau\frac{d}{d\tau}$:
 if we write $D=\tau\frac{d}{d\tau}+\sum_{n\ge 0}\tau^n D_n$, where $D_n$ is a vector field on the smooth part of $C$, then the covariant derivative is induced by $T_\gfrak(D)=\tau\frac{d}{d\tau}+
\sum_{n\ge 0}\tau^n T_\gfrak(D_n)$ acting on $\HH_\ell (V)[[\tau]]=\Hcal_\ell (V\!/\!\Delta)_{\Ccal}$. From the last clause of Lemma \ref{lemma:invariant} we get
that when $U\in \HH_\ell (V)[[\tau]]$,  
\begin{multline*}
T_\gfrak(D) E_\lambda(U)=T_\gfrak(D)(U\varepsilon^\lambda)=\\
=T_\gfrak(D)(U) \varepsilon^\lambda  -\frac{c(\lambda)}{2(\ell+\check{h})}U \varepsilon^\lambda  
=E_\lambda T_\gfrak(D)(U)-\frac{c(\lambda)}{2(\ell+\check{h})}E_\lambda(U).
\end{multline*}
Since $T_\gfrak(D)$ acts on $\HH_\ell (\tilde V_{\lambda,\lambda^*})_{\tilde C}[[\tau]]$ 
simply as  derivation by $\tau\frac{d}{d\tau}$, the last clause follows.
\end{proof}

\begin{corollary}\label{cor:doublept}
The monodromy of the WZW connection acting on $\Hcal_\ell (V\!/\!\Delta)_{\Ccal}$
has finite order and acts in the summand 
$\HH _\ell (\tilde V_{\lambda,\lambda^*})_{\tilde C}$ as 
multiplication by $\exp(\pi\sqrt{-1}\frac{c(\lambda)}{\ell+\check{h}})$. 
\end{corollary}
\begin{proof}
The multivalued  flat sections of $\Hcal_\ell (V\!/\!\Delta)_{\Ccal}$
decompose under $E$ as a direct sum labeled by $P_\ell$. The summand
corresponding to $P_\ell$ is the set of solutions of the differential equation
$\tau\frac{d}{d\tau}U+\frac{c(\lambda)}{2(\ell+\check{h})}U=0$. These are clearly
of the form $u\tau^{-c(\lambda)/2(\ell+\check{h})}$  with $u\in 
\HH _\ell (\tilde V_{\lambda,\lambda^*})_{\tilde C}$. If we let $\tau$ run over
the unit circle, then we see that the monodromy is as asserted. Finally, we observe that
$\frac{c(\lambda)}{\ell+\check{h}}\in\QQ$.
\end{proof}

\subsection*{Verlinde formula}
Theorem \ref{thm:doublept} and its Corollary \ref{cor:doublept} extend to 
the case where the base $S$ is an arbitrary $k$-variety as in Proposition
\ref{prop:fact} and the smoothing is arbitrary. This is based on a versality argument, 
which shows that our smoothing construction is not so special as it may appear.
To be concrete, suppose that we are given a family of smooth curves $\tilde\Ccal\to S$
with pairwise disjoint sections $\{S_i\}_{i\in I}\cup \{S_-,S_+\}$ and a generator
$t_{\pm}$ of the ideal defining $S_{\pm}$ in the formal completion along
$\Ccal$. Assume that the complement of the union of these sections is affine over $S$
and that this family is \emph{versal} as a family of pointed curves. 
Then $\tilde\Ccal\to S$ factors through a family whose fibers have a single double point: 
$\tilde\Ccal\to\Ccal\to S$, where  $\tilde\Ccal$ is obtained by identifying 
the sections $S_+$ and $S_-$. We regard the latter as endowed with the sections
$\{S_i\}_{i\in I}$ so that $(\Ccal,\{S_i\}_{i\in I})\to S$ is a family of pointed curves.
Then the smoothing of $\Ccal$ over $S\times\Delta$ defined by $t_\pm$ with its sections
$\{S_i\times\Delta\}_{i\in I}$ will be a versal (as a family of pointed curves) so 
that any  deformation of $(\Ccal,\{S_i\}_{i\in I})\to S$ is obtained from this one 
by means of a base change. As any two versal deformations are isomorphic, it follows that 
\ref{thm:doublept} and \ref{cor:doublept} apply to any versal deformation of 
$(\Ccal,\{S_i\}_{i\in I})\to S$.

The preceding leads to a formula of the dimension of a 
conformal block. By theorem  \ref{thm:doublept}, or rather the generalization 
discussed above, the  dimension of a conformal block stays the same under a 
degeneration. Since every pointed curve degenerates into one with ordinary 
double points whose normalization consists of curves of genus zero, it suffices 
to do the computation for such a degenerate curve. But then we may invoke
\ref{prop:fact} to reduce to the case of a smooth rational curve, which can
be dealt with using our discussion of the genus zero case. We can even arrange that
the support of the representation map meets every component in at most three points.
A more refined approach involves the notion of a fusion ring \cite{beauville}.

\begin{remark}
In case $k=\CC$, one can work in the complex-analytic category. Then
the fact that the singularity is formally regular singular ensures
that the flat multivalued sections converge on simply connected sectors
based at $o$ (see for instance \cite{malgrange}) so that the monodromy has 
the expected interpretation. If $C-\{S_0\}$ is smooth, then the 
we have an associated fibration over the punctured unit disk whose monodromy
is given by a Dehn twist along a vanishing circle on the general fiber.
In terms of the setting of Remark \ref{rem:topmeaning}: 
an isotopy class $\alpha$ of embedded circles in $\Sigma$ defines the class
of a Dehn twist $D_\alpha\in\Gamma(\Sigma ,I)$ and the action of $D_\alpha$ on
$\PP (\Sigma,I,V)$ is given by the above formula (hence is of finite order).
\end{remark}

\end{document}